\documentclass[a4paper,12pt]{article}
\usepackage[american]{babel}
\usepackage{amsmath}
\usepackage{latexsym}
\usepackage{amssymb}
\usepackage{theorem}
\usepackage{diagrams}
\theoremstyle{change}
\newtheorem{Thm}{Theorem}[section]
\newtheorem{Cor}[Thm]{Corollary}
\newtheorem{Prop}[Thm]{Proposition}
\newtheorem{Lem}[Thm]{Lemma}
{\theorembodyfont{\rmfamily}

}
\renewcommand{\phi}{\varphi}

\newcommand{\proof}{\par\medskip\rm\emph{Proof. }}

\newcommand{\qed}{\ \hglue 0pt plus 1filll $\Box$}

\newcommand{\RR}{\mathbb{R}}
\newcommand{\ZZ}{\mathbb{Z}}

\renewcommand{\SS}{\mathbb{S}}

\newcommand{\SKIP}[1]{}
\newcommand{\eps}{\varepsilon}
\renewcommand{\emptyset}{\varnothing}

\begin{document}

\title{{\bf Completely reducible subcomplexes of spherical buildings}}
\author{Linus Kramer}
\date{April 15 2008}
\maketitle

\begin{abstract}
A completely reducible subcomplex of a spherical building is a spherical building.
\end{abstract}
By a \emph{sphere} we mean metric space isometric to the unit sphere
$\SS^m\subseteq\RR^{m+1}$, endowed with the spherical
metric $d$. The distance of $u,v\in\SS^m\subseteq\RR^{m+1}$ is given
by $\cos(d(u,v))=u\cdot v$ (standard inner product).
Recall that a geodesic metric space $Z$ is CAT$(1)$ if no geodesic
triangle of perimeter $<2\pi$ in $Z$ is thicker than its
comparison triangle in the $2$-sphere $\SS^2$.  It follows that
any two points at distance $<\pi$ can be joined by a unique geodesic segment
(an isometric copy of a closed interval).
A subset $Y$ of a CAT(1)-space $Z$ is \emph{convex} if the following
holds: for all $x,y\in Y$ with $d(x,y)<\pi$, the geodesic segment
$[x,y]$ is contained in $Y$. It is clear from the definition
that arbitrary intersections of convex subsets are convex
(and CAT$(1)$).

\section{Convex sets in Coxeter complexes}
Let $\Sigma$ be an $n$-dimensional spherical Coxeter complex, let
$\bar\Sigma$ be a simplicial complex which
refines the triangulation of $\Sigma$ and which is
invariant under the Coxeter group  $W$ and $\pm id$. Examples of such triangulations are
$\Sigma$ itself and its barycentric subdivisions. In the geometric
realization, the simplices are assumed to be spherical.
The \emph{span} of a subset of a sphere is the smallest subsphere
containing the set.

We assume now that $A\subseteq\bar\Sigma$ is an $m$-dimensional
subcomplex whose geometric realization $|A|$ is convex.
\begin{Lem}
Let $a\in A$ be an $m$-simplex. Then
$|A|\subseteq \mathrm{span}|a|$.
\proof
Assume this is false. Let $u\in|A|\setminus\mathrm{span}|a|$.
Then $-u\not\in|a|$ and
$Y=\bigcup\{[u,v]\mid v\in|a|\}$ is contained in $|A|$. But
$Y$ is a cone over $|a|$ and in particular $m+1$-dimensional,
a contradiction.
\qed
\end{Lem}
We choose an $m$-simplex $a\in A$ and put 
\[
S=\mathrm{span}|a|\cap|\bar\Sigma|;
\]
this is an $m$-sphere containing $|A|$.
Recall that an $m$-dimensional simplicial complex is called
\emph{pure} if every simplex is contained in some $m$-simplex.
\begin{Lem}
$A$ is pure.

\proof
It suffices to consider the case $m\geq 1$.
Let $a\in A$ be an $m$-simplex, and assume that $b\in A$ is a
lonely simplex of maximal dimension $\ell<m$.
Then $int(-b)$ is disjoint from $int(a)$. Let
$v$ be an interior point of $a$ and $u$ an interior point
of $b$ and consider the geodesic segment
$[u,v]\subseteq|A|$.
If $b$ is a point, the existence of the geodesic
shows that $b$ is contained in some higher dimensional simplex,
a contradiction. If $\ell\geq 1$, then $[u,v]$ intersects
$int(b)$ in more than two points (because $b$ is lonely),
so $v$ is in the span of $b$. This contradicts
$dim(a)>dim(b)$.
\qed
\end{Lem}
\begin{Lem}
\label{OppLem}
If there exists an $m$-simplex $a\in A$ with $-a\in A$, then
$|A|=S$.

\proof
Then any point in $S$ lies on some geodesic of length $<\pi$
joining a point in $|a|$ with a point in $|-a|$.
\qed
\end{Lem}
Topologically, the convex set
$|A|$ is either an $m$-sphere or homeomorphic
to a closed $m$-ball. For $m\geq 2$, these spaces are strongly
connected (i.e. they cannot be separated by $m-2$-dimensional
subcomplexes \cite{Alex}). It follows that $A$ is a chamber complex,
i.e. the chamber graph $C(A)$ (whose vertices are the $m$-simplices
and whose edges are the $m-1$ simplices) is connected
\cite{Alex}. If $m=1$, then $|A|$ is a connected graph and
hence strongly connected.
\begin{Lem}
If $m\geq 1$, then $A$ is a chamber complex.
\qed
\end{Lem}

\section{Results by Balser-Lytchak and Serre}
We now assume that $X$ is a simplicial spherical building
modeled on the Coxeter complex $\Sigma$. By means of the coordinate
charts for the apartments we obtain a metric simplicial complex
$\bar X$ refining $X$, which is modeled locally on
$\bar\Sigma$. In this refined complex $\bar X$, we call two simplices
$a,b$ \emph{opposite} if $a=-b$ in some (whence any) apartment
containing both. We let $opp(a)$ denote the collection of all
simplices in $\bar X$ opposite $a$.
The geometric realization $|\bar X|$ is CAT$(1)$.
Furthermore, any geodesic arc is contained in some apartment.

We assume that $A\subseteq\bar X$ is
an $m$-dimensional subcomplex and that $|A|$ is convex.
For any two simplices $a,b\in A$, we can find an apartment $\bar\Sigma$
containing $a$ and $b$. The intersection $|A|\cap|\bar\Sigma|$ is then
convex, so we may apply the results of the previous section
to it. We note also that $|A|$ is CAT$(1)$.
\begin{Lem}
$A$ is a pure chamber complex.

\proof
Let $a\in A$ be an $m$-simplex and let $b\in A$ be any simplex.
Let $\bar\Sigma$ be an apartment containing $a$ and $b$.
Since $|\bar\Sigma|\cap|A|$ is $m$-dimensional and
convex, we find an $m$-simplex $c\in A\cap\bar\Sigma$
containing $b$. Similarly we see that $A$ is a chamber complex.
\end{Lem}
The next results are due to Serre \cite{Serre} and Balser-Lytchak \cite{BL1,BL2}.
\begin{Lem}
If there is a simplex $a\in A$ with $opp(a)\cap A=\emptyset$,
then $|A|$ is contractible.

\proof
We choose $u$ in the interior of $a$. Then
$d(u,v)<\pi$ for all $v\in|A|$, so $|A|$ can be contracted to
$u$ along these unique geodesics.
\qed
\end{Lem}
\begin{Prop}
\label{SerreProp}
If there is an $m$-simplex $a$ in $A$ with $opp(a)\cap A\neq\emptyset$,
then every simplex $a\in A$ has an opposite in $A$.

\proof
Let $a,b\in A$ be opposite $m$-simplices, let $\bar\Sigma$ be
an apartment containing both and let $S\subseteq|\bar\Sigma|$ denote the
sphere spanned by $a,b$. Then $S\subseteq|A|$.
Let $c$ be any $m$-simplex in $A$.
If $c$ is not opposite $a$, we find interior points
$u,v$ of $c,a$ with $d(u,v)<\pi$.
The geodesic arc $[u,v]$ has a unique extension in
$S$. Along this extension, let $w$ be the point with $d(u,w)=\pi$ and
let $c'$ be the smallest simplex containing $w$. Then $c'$ is opposite $c$.

Thus every $m$-simplex in $A$ has an opposite, and therefore every
simplex in $A$ has an opposite.
\qed
\end{Prop}
In this situation where every simplex has an opposite, $A$ is called
$A$ \emph{completely reducible}.
If every simplex of a fixed dimension $k\leq m$ has an opposite in $A$,
then clearly every vertex in $A$ has an opposite. Serre \cite{Serre}
observed that the latter already characterizes complete reducibility.
\begin{Prop}
\label{SerreProp2}
If every vertex in $A$ has an opposite, then $A$ is completely reducible.

\proof
We show inductively that $A$ contains  a pair of opposite $k$-simplices,
for $0\leq k\leq m$.
This holds for $k=0$ by assumption, and we are done if $k=m$ by
\ref{SerreProp}. So we assume that $0\leq k<m$.

Let $a,a'$ be opposite $k$-simplices in $A$ and let $b\in A$ be a vertex
which generates together with $a$ a $k+1$-simplex
(recall that $A$ is pure, so such a vertex exists).
We fix an apartment $\bar\Sigma$ containing $a$, $b$ and $a'$.
The geodesic convex closure $Y$ of $b$ and $|a|\cup|a'|$ in the sphere
$|\bar\Sigma|$ is a $k+1$-dimensional hemisphere (and is
contained in $|A|$). Let $b'\in A$ be a vertex opposite $b$.
A small $\eps$-ball in $Y$ about $b$ generates together
with $b'$ a $k+1$-sphere $S\subseteq|A|$. Because
$\dim S=k+1$, there exists a point $u\in S$ such that the minimal
simplex $c$ containing $u$ has dimension at least $k+1$.
Let $u'$ be the opposite of $u$ in $S$, and $c'$ the minimal simplex
containing $u'$. Then $c,c'$ is a pair of opposite simplices in $A$
of dimensions at least $k+1$.
\qed
\end{Prop}

\section{Completely reducible subcomplexes are buildings}
We assume that $A$ is $m$-dimensional, convex and completely
reducible. If $m=0$, then $A$ consists of a set of
vertices which have pairwise distance $\pi$. This set is, trivially,
a $0$-dimensional spherical building. So we assume now
that $1\leq m\leq n$.
Two opposite $m$-simplices $a,b\in A$ determine an $m$-sphere 
$S(a,b)$ which we call a \emph{Levi sphere}.
\begin{Lem}
If $a,b\in A$ are $m$-simplices, then there
is a Levi sphere containing $a$ and $b$.

\proof
This is true if $b$ is opposite $a$. If $b$ is not
opposite $a$, we choose interior points $u\in int(a)$ and
$v\in int(b)$, and  a simplex $c\in A$ opposite $b$.
The geodesic $[u,v]$ has a unique continuation 
$[v,w]$ in the Levi sphere $S(b,c)$, such that $d(u,w)=\pi$.
Let $\bar\Sigma$ be an apartment containing
the geodesic arc $[u,v]\cup[v,w]$ and let $d$ be the smallest simplex
in $\bar\Sigma$ containing $w$. Then $d$ is in $A$ and opposite
$a$, so the there is a
Levi sphere $S(a,d)$ containing $[u,v]\cup[v,w]$. Since $b$ is
the smallest simplex containing $v$, it follows 
that $b\in S(a,d)$.
\qed

\end{Lem}
Since $A$ is pure, we have the following consequence.
\begin{Cor}
Any two simplices $a,b\in A$ are in some Levi sphere.
\qed
\end{Cor}
We call an $m-1$-simplex $b\in A$ \emph{singular} if it is contained
in three different $m$-simplices. The following idea is taken from
Caprace \cite{Cap}. Two $m$-simplices are \emph{t-equivalent}
if there is a path between them in the dual graph 
which never crosses a singular $m-1$-simplex. The t-class of
$a$ is contained in all Levi spheres containing $a$.
\begin{Lem}
Let $b$ be a singular $m-1$-simplex. Let $S$ be a Levi sphere
containing $b$ and let $H\subseteq S$ denote the great
$m-1$-sphere spanned by $|b|$.
Then $H$ is the union of singular $m-1$-simplices.

\proof
Let $a$ be an $m$-simplex containing $b$ which is not in $S$ and
let $-b$ denote the opposite of $b$ in $S$. Let $S'$ be
a Levi sphere containing $a$ and $-b$ and consider the
convex hull $Y$ of $|a|\cup|-b|$ in $S'$.
Then $Y$ is an $m$-hemisphere. The intersection $Y\cap S$ is convex,
contains the great sphere $H$, and is different from $Y$,
so $Y\cap S=H$.
\qed
\end{Lem}
We call $H$ a \emph{singular great sphere}. Along singular great spheres,
we can do 'surgery':
\begin{Lem}
Let $S,H,Y$ be as in the previous lemma. Let $Z\subseteq S$ be a
hemisphere with boundary $H$. Then $Z\cup Y$ is a Levi sphere.

\proof
We use the same notation as in the previous lemma. Let
$c\subseteq Z$ be an $m$-simplex containing $-b$, then
$|c|\cup H$ generates $Z$. Let $S'$ be a Levi sphere containing
$c$ and $a$. Then
$Z\cup Y\subseteq S'$ and $Z\cup Y=H$, whence
$S'=Z\cup Y$.
\qed
\end{Lem}
\begin{Lem}
Let $S$ be a Levi sphere and let $H,H'\subseteq S$ be singular
great spheres. Let $s$ denote the metric reflection of
$S$ along $H$. Then $s(H')$ is again a singular great sphere.

\proof
We use the notation of the previous lemma.
Let $b'$ be a singular $m-1$-simplex in $H'\cap Z$. Let
$-b'$ denote its opposite in the Levi sphere $S'=Z\cup Y$.
We note that the interior of $b$ is disjoint from $S$.
Let $b''$ be the opposite of $-b'$ in the Levi sphere
$S''=(S\setminus Z)\cup Y$. Then $b''$ is a singular
$m-1$-simplex in $S$, and $b''$  is precisely the
reflection $s(b')$ of $b'$ along $H$.
\qed
\end{Lem}
For every Levi sphere $S$ we obtain in this way
a finite reflection
group $W_S$ which permutes the singular great spheres in
$S$. As a reprentation sphere, $S$ may split off a trivial
factor $S_0$, the intersection of all
singular great spheres in $S$. We let $S_+$ denote its orthogonal
complement, $S=S_0*S_+$.
The intersections of the singular great spheres with $S_+$
turn $S_+$ into a spherical Coxeter complex, with Coxeter
group $W_S$. Let $F\subseteq S$ be a fundamental domain for
$W_S$, i.e. $F=C*S_0$, where $C\subseteq S_+$ is a Weyl chamber.
The geometric realization of the t-class of any $m$-simplex
in $F$ is precisely $F$.
\begin{Lem}
If two Levi spheres $S,S'$ have an $m$-simplex $a$ in common, then there
is a unique isometry $\phi:S\rTo S'$ fixing $S\cap S'$ pointwise.
The isometry fixes $S_0$ and maps $W_S$ isomorphically onto $W_{S'}$.

\proof
The intersection $Y=S\cap S'$ contains the fundamental domain
$F$. Since $F$ is relatively open in $S$, there
is a unique isometry $\phi:S\rTo S'$ fixing $Y$. The Coxeter group
$W_s$ is generated by the reflections along the singular $m-1$-simplices
in $Y$. Therefore $\phi$ conjugates $W_S$ onto $W_{S'}$.
Finally, $Y$ contains $S_0$.
\qed
\end{Lem}
\begin{Cor}
If two Levi spheres $S,S'$ have a point $u$ in common, then
$S_0=S_0'$. Furtheremore, there exists an isometry $\phi:S\rTo S'$
which fixes $S\cap S'$ and which conjugates $W_S$ to $W_{S'}$.

\proof
Let $a,a'$ be $m$-simplices in $S$ and $S'$ containing $u$, and let
$S''$ be a Levi sphere containing $a$ and $a'$.
We compose $S\rTo S''\rTo S'$.
\end{Cor}
\begin{Thm}
Let $A$ be completely reducible. Then there is a thick spherical building
$Z$ such that $|A|$ is the metric realization of $Z*\SS^0*\cdots*\SS^0$.

\proof
Let $S$ be a Levi sphere and let $k=\dim S_0+1$. We make
$S_0$ into a Coxeter complex with Coxeter group $W_0=\ZZ/2^k$
(we fix an action, this is not canonical). By the previous
Corollary, we can transport the simplicial structure on
$S$ unambiguously to any Levi sphere in $A$.
\qed
\end{Thm}
For $A=X$, this is Scharlau's reduction theorem for weak
spherical buildings \cite{Scha} \cite{Cap}.

\end{document}